\documentclass{fundam}

\usepackage[T2A,T1]{fontenc}
\usepackage[utf8]{inputenc}
\usepackage[russian,english]{babel}

\usepackage{graphicx}
\usepackage{amssymb}
\usepackage[font=small]{caption}
\usepackage{floatrow}
\usepackage{setspace}
\usepackage{url}

\newcommand{\Qn}{\medskip\noindent\textbf{Q:\ }}
\newcommand{\Ar}{\medskip\noindent\textbf{A:\ }}
\newcommand{\R}{\mathbb R}
\newcommand{\Th}{\text{Th}}

\begin{document}
	
\setcounter{page}{133}
\publyear{22}
\papernumber{2122}
\volume{186}
\issue{1-4}

    \finalVersionForARXIV
    %%%  \finalVersionForIOS

\title{The 1966 International Congress of Mathematicians: \\ A Micro-memoir}

\author{Yuri Gurevich\thanks{Address for correspondence: Computer Science and Engineering, University of Michigan, USA}
             \\
Computer Science and Engineering\\
University of Michigan, USA\\
gurevich@umich.edu}

\maketitle

\runninghead{Y. Gurevich}{The 1966 International Congress of Mathematicians: A Micro-memoir}

\vspace*{-5mm}
\begin{quote}\raggedleft\noindent
To the memory of Boris (Boaz) Trakhtenbrot, \\
a founding father of theoretical computer science\\
and a friend who is greatly missed
\end{quote}
%\bigskip

\begin{figure}[!h]
\vspace*{3mm}
\includegraphics[width=8.6cm]{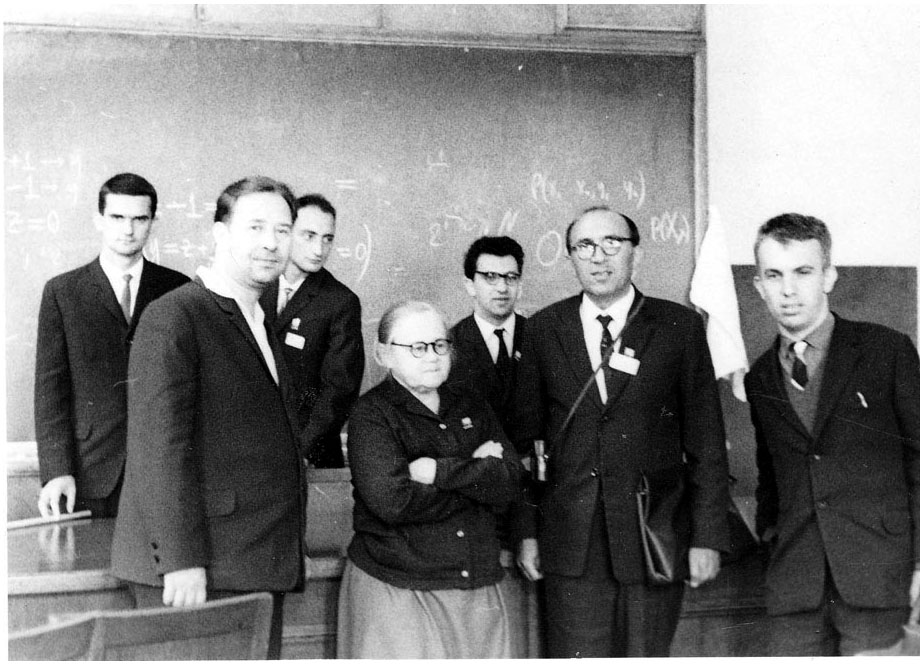}\vspace*{-1mm}
\caption{Boris Trakhtenbrot with Russian speaking colleagues at the 1966 congress}
\label{Trakh}
\end{figure}

\section{A hole in the iron curtain}
\label{sec:curtain}

\noindent\textbf{Quisani\footnotemark:}
I guess an International Congress of Mathematicians is a big deal?

\footnotetext{Author's former student.}

\medskip\noindent\textbf{Author:}
It has an illustrious history \cite{ICM}. The first congress was held in Zurich in 1897. Hilbert's famous list of problems was presented at the 1900 congress, in Paris.
The 1966 congress was special for Soviet mathematicians because it was held in Moscow. We had been deprived of contacts with our Western colleagues.

\Qn Were Western books and journals available in the USSR?

\Ar There were good mathematical libraries in Moscow, Leningrad and Novosibirsk. I know of no other place with a good mathematical library. To buy Western books or to subscribe to Western journals, you needed foreign currency. The Ural State University, where I was teaching at the time, had zero foreign currency for that purpose, even though it was the main university of the whole area of the Urals.

\Qn Could you move to Moscow?

\Ar For all practical purposes, no. The Soviet system was feudal in some respects. The right to live in Moscow normally went from parents to their children. Even if some Moscow institution was eager to hire you, in most cases they would not be able to do that \cite{Propiska}.

\Qn Could you talk to foreigners?

\Ar That was risky. Besides, there weren't  many foreigners to talk to. The Urals were completely closed to foreigners.

\Qn To all foreigners or only to Westerners?

\Ar I believe the Urals were closed to all foreigners, though there could have been some tightly controlled exceptions. I remember that, right after the congress, the Polish logician Bogdan W\c eglorz and I considered flying together to Novosibirsk, but I needed to stop in Sverdlovsk (the cultural center of the Urals, now called Yekaterinburg), and Bogdan was not allowed to stop there.

\Qn At the congress, were your Western colleagues willing to talk to you?

\Ar Yes. They --- and especially the Americans --- were very smiley. Of course I smiled back, but my smiling muscles were not sufficiently developed. At the end of the congress, those muscles were aching.

\Qn Were you particularly gloomy?

\Ar I don't think so. We, in Russia, just didn't smile that often, certainly not at strangers. ``A laugh without a reason is a sign of folly,'' says a Russian adage.%
\footnote{\selectlanguage{russian}%
Смех без причины --- признак дурачины.}

\Ar Did you have time to speak to your Soviet colleagues?

\Qn Of course I spoke to my Soviet colleagues, though mostly to my peers.
On one occasion, running from one building to another, I met Anatoly I. Malcev, the father of Soviet model theory. Malcev was super polite, and we had a small talk about the difficulty of disseminating one's own results. ``Find a famous guy,'' joked Malcev, ``grab him by a button and force feed him your results.''

\Qn Did you know Trakhtenbrot at the time?

\Ar Yes, I knew him, but I do not remember talking to him at the congress.
I had an opportunity to attend his Novosibirsk seminar a couple of times.
It was a pioneering seminar on nascent computer science.
But, at the time, I was moving from algebra to logic, and I did not realize how forward-looking Trakhtenbrot's seminar was.

\section{Logic}
\label{sec:logic}

\Qn The congress, I am sure, was very interesting to you scientifically. Give me a couple of highlights.

\begin{figure}[!h]
\vspace{4mm}
\includegraphics[width=7.8cm]{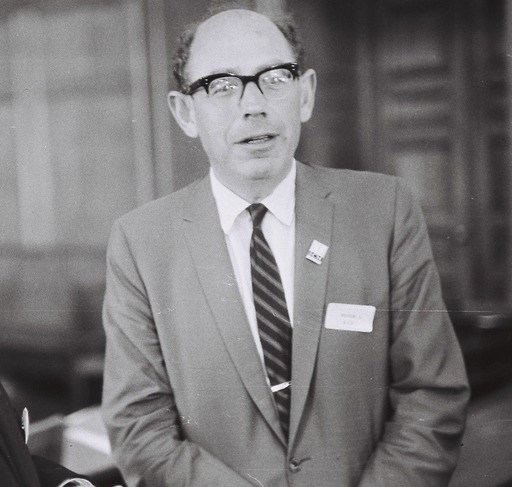}
\caption{Abraham Robinson at the congress}\label{RobinsonPic}
\end{figure}

\Ar One highlight for me was a beautiful lecture on nonstandard analysis given by Abraham Robinson, the originator of that field. Infinitesimals used to be a fundamental notion in mathematical analysis. Eventually, for the purpose of rigor, their use was replaced by $\varepsilon$-$\delta$ arguments. Intuitively a positive infinitesimal $\iota$ is an infinitely small real number different from zero, so that we have
\begin{equation}\label{iota}
0<\iota<1/n\quad \text{for all positive integers}\ n.
\end{equation}

\Qn But this makes no sense.

\Ar Indeed, it seems so. Yet nonstandard analysis provides a solid foundation for (a careful use of) infinitesimals.

\Qn Explain.

\Ar Extend the first-order theory $\Th(\R)$ of the ordered real field $\R$ with a constant symbol $\iota$ and the axioms \eqref{iota}. Every finite subset of these axioms is obviously consistent with $\Th(\R)$. By the compactness theorem, the set of all these axioms is consistent with $\Th(\R)$ and thus has a model $\R^*$. In that model, the value of $\iota$ is a positive nonstandard real number that is less than any standard positive real number.

With infinitesimals, one can make rigorous the intuitive ideas underlying analysis. For example, a real-valued function $f$ of a real variable is continuous at a real number $a$ if and only if an infinitesimal change of the input produces an infinitesimal change of the output, i.e., for every infinitesimal $\xi$, the difference $f(a+\xi) - f(a)$ is infinitesimal. Moreover, if the ratio $\displaystyle\frac{f(a+\xi) - f(a)}\xi$ is infinitely close to a standard real number $b$ for all nonzero infinitesimals $\xi$ then $b$ is the derivative of $f$ at $a$.

At the end of his lecture, Robinson announced that he would give away a few copies of his book \cite{Robinson}. I approached him and asked whether I can have one. Only if you promise, he said, that you will read the book. I told him that I would do better, that I would conduct a seminar on nonstandard analysis. He gave me a copy of his book, and I indeed conducted such a seminar at the Ural State University.

\Qn Was there anything at the congress that impressed you even more than nonstandard analysis?

\begin{figure}[!h]
\vspace*{3mm}
\includegraphics[width=7.8cm]{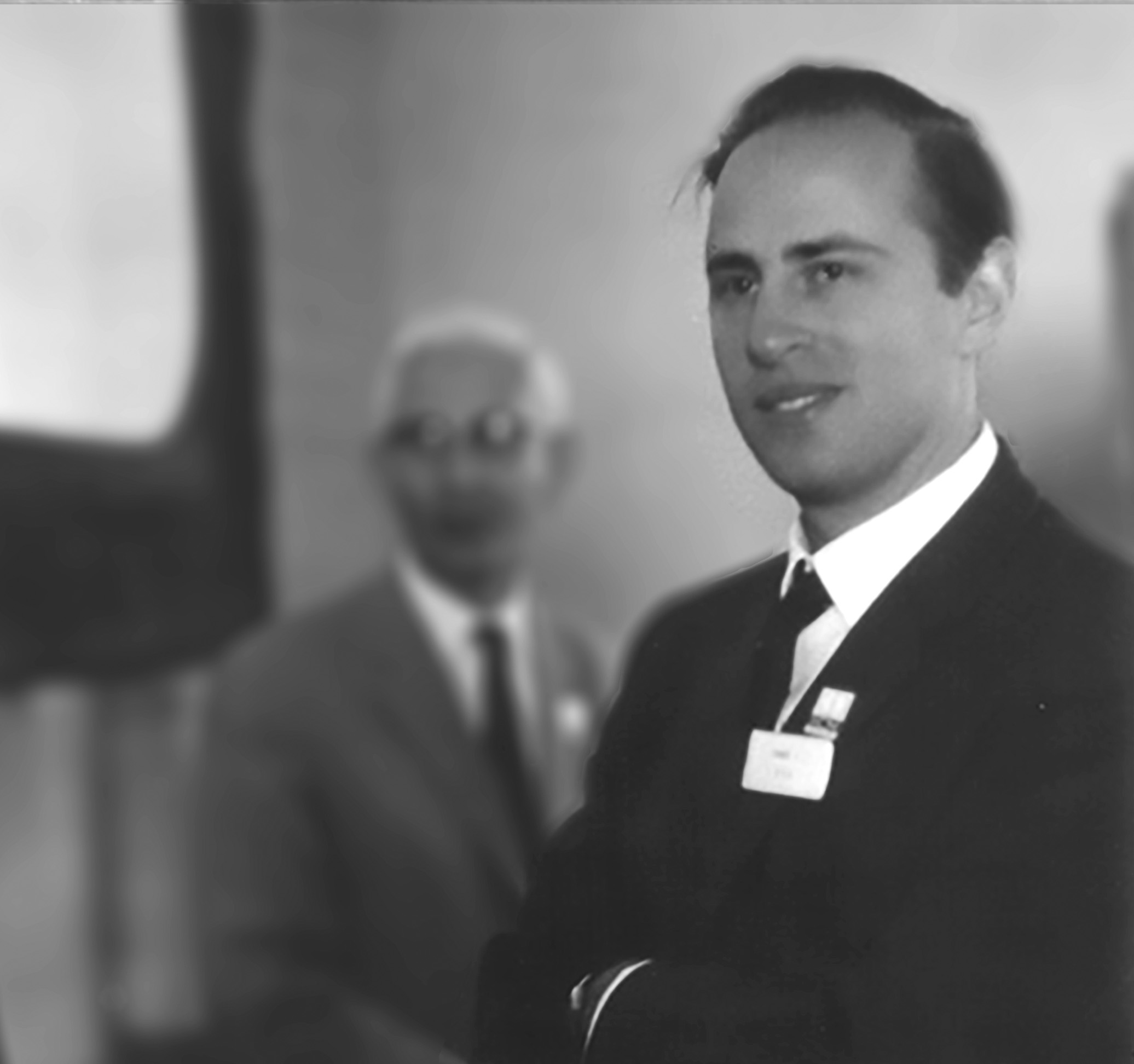}
\caption{Paul Cohen at the congress}\label{CohenPic}
\end{figure}

\Ar Yes. This was the method of forcing presented at the congress by its inventor, Paul Cohen. Forcing allows you, given a set-theoretic world $W$, to construct new set-theoretic worlds by adding elements to $W$ in a tightly controlled way.

\Qn What is a set-theoretic world?

\Ar Typically, it is a model of ZFC, the first-order Zermelo-Fraenkel set theory with the axiom of choice.

\Qn This seems similar to nonstandard analysis. There you extend a given model $\R$ of the first-order theory $\Th(\R)$ while here you extend a given model $W$ of the first-order theory ZFC. I guess the difference is that ZFC is more complicated than $\Th(\R)$.

\Ar There is something else. You get a lot of mileage from an arbitrary non-standard model $\R^*$ of $\Th(\R)$. Any such $\R^*$ allows you to formalize nicely many intuitively appealing arguments that involve infinitesimals. In the set-theoretic case, an arbitrary extension of the given model doesn't buy you much. You need models with particular properties. Cohen used forcing to construct a model of ZFC where the continuum hypothesis (CH) fails \cite{CohenBook}. Earlier, G\"odel had constructed a model of ZFC + CH \cite{Goedel}. Thus Cohen proved the independence of the continuum hypothesis.

The independence result all by itself was an enormous achievement; at the congress, Cohen received a Fields medal \cite{Fields} and in general was a big hero. But the forcing method was more important yet. Very quickly it radically changed set theory. A great many problems in set theory, algebra, topology, etc.\ have been solved using forcing.

\Qn Do you understand the method of forcing?

\Ar I learned forcing in the 1970s when I attended the logic seminar at the Hebrew University of Jerusalem. It was virtually impossible to be an active member of the seminar without understanding forcing.

But, during Cohen's lecture at the congress, I did not understand a thing. After his lecture, a group of Soviet mathematicians, about two dozen of us, asked Cohen to give us an informal lecture. Cohen  spoke, in English, and took numerous questions. One of the Moscow mathematicians was translating. We stopped when, after a couple of hours, the translator was exhausted. We did not have a substitute translator from English. Cohen knew some French and German, but we did not have translators from those languages either.

As a result of that informal session I understood better why ZFC + CH was consistent. I had tried to study G\"odel's book on the subject but it was too formal and dry. Cohen's explanation made G\"odel's approach clear. But I did not understand why ZFC + the negation of CH was consistent. The whole idea of forcing looked mystical.

\Qn Maybe I should read Cohen's book \cite{CohenBook}.

\Ar This will not be easy. Cohen gives a simple and beautiful proof of G\"odel's result, but the proof of his own result is challenging. A streamlined exposition of Cohen's result was given by Shoenfield \cite{Shoenfield}.

\Qn Did you speak to Cohen?

\Ar I wouldn't dare to but he approached me and asked whether I am Yuri Gurevich. I said yes, I am, but I do not know anything about forcing. He smiled and said that he had been working on quantifier elimination for some field theories and that he wanted to know how I proved the decidability of the theory of ordered Abelian groups. I sketched my proof.

\Qn Were there many logicians at the congress?

\Ar Everybody was there, it seems, except for Kurt G\"odel. The authors of the logic books used to be just names for me. But there they were: Alonzo Church, Stephen Kleene, Alfred Tarski.

\begin{figure}[ht]
\vspace*{3mm}
\includegraphics[width=9cm]{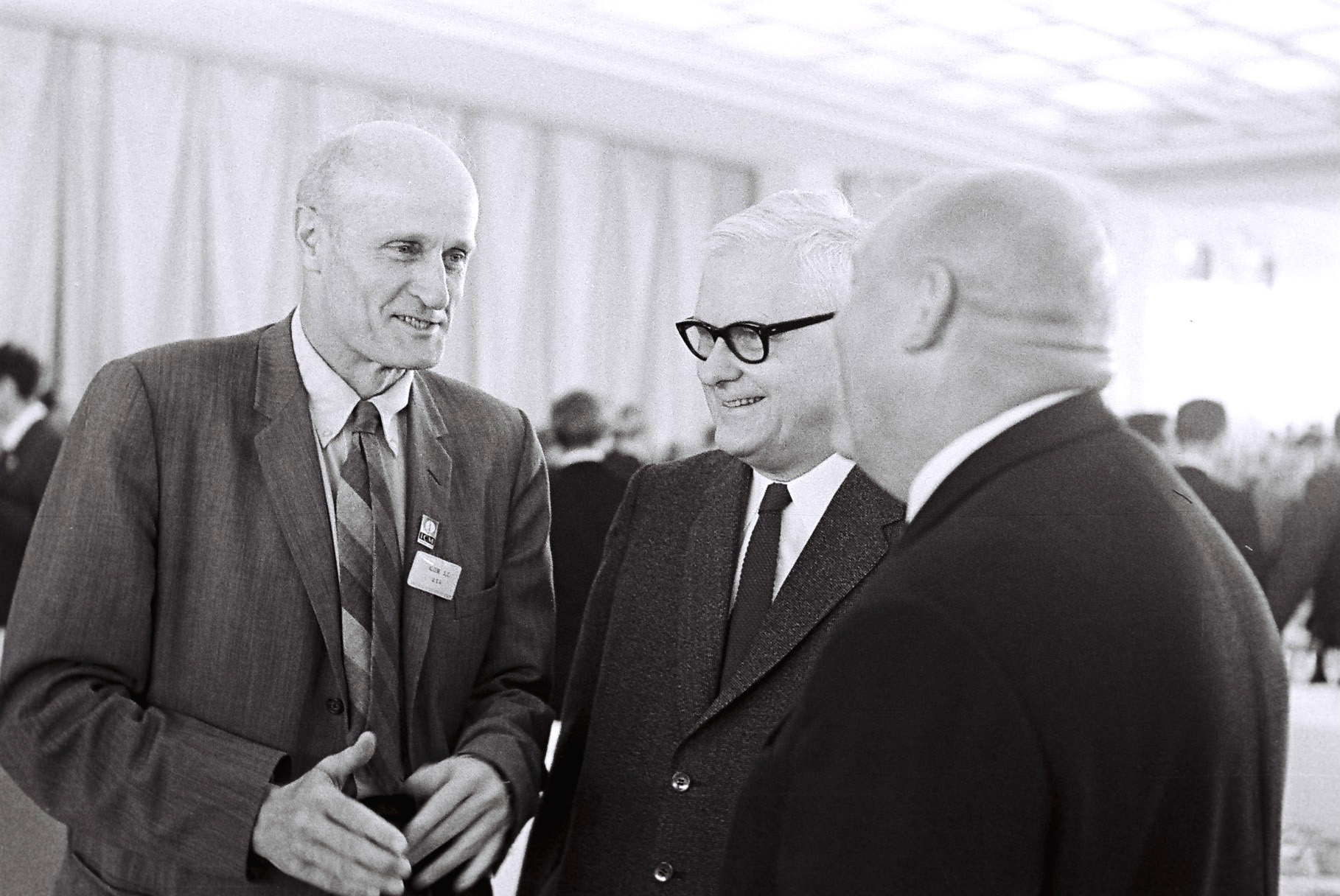}
\caption{Stephen Kleene, Alonzo Church and Anatoly Malcev at the congress}\label{Kleene}
\end{figure}

\begin{figure}[h]
\vspace*{3mm}
\includegraphics[width=9cm]{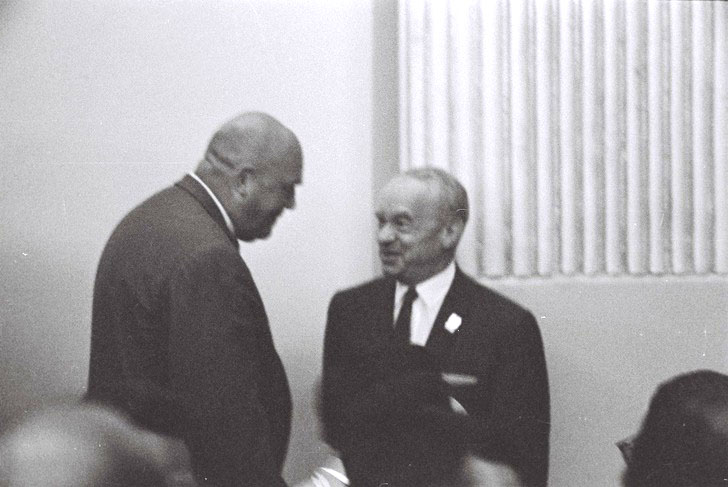}
\caption{Anatoly Malcev and Alfred Tarski at the congress}\label{Tarski}
\end{figure}

\Qn Did you give a talk at the congress?

\Ar I gave a 15 minute talk, on the classical decision problem \cite{Gurevich}. Translators had not been provided for short presentations but Haim Gaifman, of the Hebrew University of Jerusalem, helped me out. I split the talk into several pieces and explained every piece to Haim beforehand. During the presentation, every piece was presented first in Russian by me and then --- loud and clear --- in English by Gaifman.

\begin{figure}[h]
\vspace*{3mm}
\includegraphics[width=8cm]{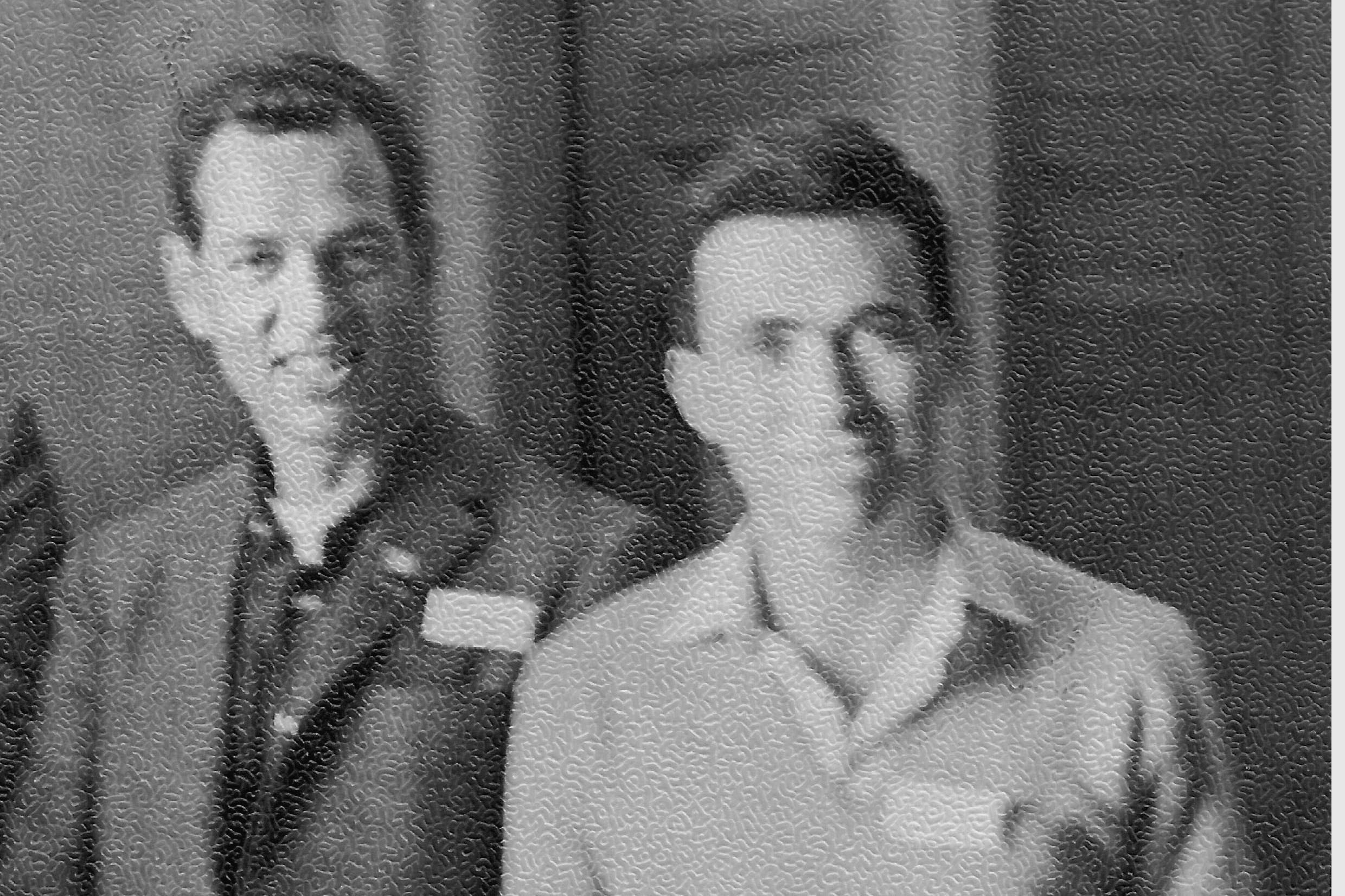}
\caption{Sol Feferman and Haim Gaifman at the congress}\label{Feferman}
\end{figure}

\section{Language}
\label{sec:lang}

\Qn How did you communicate with Robinson, Cohen and Gaifman? Did they know some Russian?

\Ar No, they didn't. Russian is an important language in mathematics, and many Western mathematicians could read Russian mathematical texts. But, in my experience, very few of them spoke Russian. One exception was Alfred Tarski who spoke flawless Russian. Warsaw, his home town, was in the Russian Empire during his first 14 years.

\Qn How was your English?

\Ar Nonexistent. I was a monolingual Russian speaker who never spoke any other language whatsoever. But burning desire is worse than fire.%
\footnote{\selectlanguage{russian}Oхота пуще неволи.}
I used an unprincipled concoction of German that I studied in school but never used and Yiddish that I heard at home. My parents spoke Yiddish between themselves but not to me. At the congress, my communication skills grew with use. On at least one occasion, I even translated a conversation between a Russian and an American.  One of my American interlocutors joked that there were three universal languages: mathematics, music and Yiddish.

One peculiar problem was related to the fact that G and H are often rendered by the same Russian letter
\selectlanguage{russian}Г.
\selectlanguage{english}``Hegel'' becomes
\selectlanguage{russian}``Гегель.''
\selectlanguage{english}%
You have to guess whether it was ``Gegel, Gehel, Hegel,'' or ``Hehel,'' At one point, I asked Paul Cohen, in a written form, whether he met ``H\"odel.'' He looked puzzled but then figured out my intent and said that he did meet G\"odel.

\Qn Why the written form?

\Ar I mentioned above my conversation with Cohen. We spoke for a while, and then Cohen said that he wanted to go to a certain lecture. I went along, and we sat together. After a minute or two, he lost interest in the lecture, and I initiated the written exchange.

\begin{figure}[hb!]
\vspace*{3mm}
\includegraphics[width=6cm]{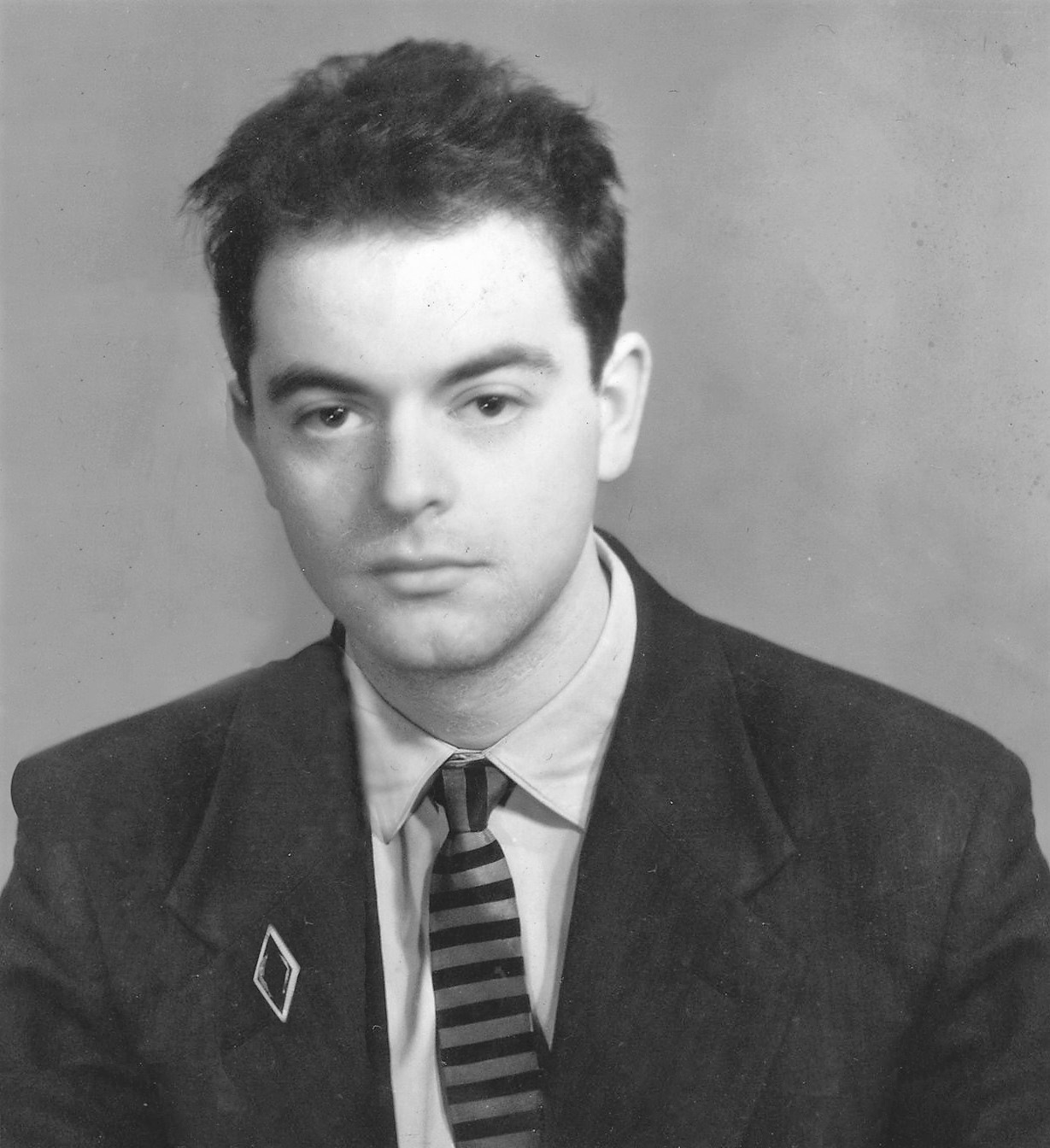}
\caption{The author shortly before the congress} \label{Yuri}
\end{figure}

\section{Epilog}
\label{sec:epilog}

\Qn Were you afraid of possible consequences of your interaction with Westerners at the congress?

\Ar I was, of course. The Soviet Union was changing unpredictably after Stalin's death in 1953. The Khrushchev Thaw was a welcome development but it was inconsistent. Besides, Khrushchev was removed from power in 1964. It wasn't clear at all in 1966 which way the country would go.

\Qn Were all foreigners equally dangerous for you to talk to?

\Ar At the time I thought that, for me, it was most dangerous to talk to Israelis. I was happy that Haim Gaifman agreed to help me with my talk; that gave me an excuse to talk to him. By the way, I asked him whether he had a car. He said yes but not right now because his car hit some object. He named that object in many languages, and I still had no idea what it was until he waved his hands: it was a camel.

\Qn After the congress, could you correspond with your Western colleagues?

\Ar The rules were unclear. You were supposed to ask for permission. Now imagine yourself being a Soviet academic bureaucrat. My letter arrives to you with a request to send it to the West. Since the rules are unclear, it is safer for you to deny the request or to postpone the decision indefinitely.

Instead of seeking permission to correspond, I went to the post office and sent my letters, hoping for the best. Amazingly, the letters usually made it through, and I received some letters from the West.

\Qn I guess that the 1966 congress remains an unforgettable event for you.

\Ar Oh yes, it was a great learning experience. Also, I made many friends. When I was leaving the USSR in 1973, Sol Feferman, whom I met at the congress, invited me to Stanford ``if you don't go to Israel.'' When I was seeking my first job and needed recommendation letters, my Soviet colleagues were unavailable. Abraham Robinson and Alfred Tarski graciously wrote for me.

\subsection*{Acknowledgement}
The pictures in Figures~1 and 6 are from the archive of Boris Trakhtenbrot. Those in Figures~2--5 were taken, at the congress, by Prof. Sergey V. Smirnov (1911-1979) of Ivanovo State Pedagogical Institute (now Ivanovo State University), and are published here with kind permission of his daughter Olga Smirnova. The university keeps Smirnov's {\it Nachlass}, and additional congress pictures can be found at~\cite{Smirnov}.

I am grateful to Andreas Blass and Victor Marek for commenting on the draft of this micro-memoir. (Victor Marek was one of the young Polish logicians whom I met at the congress.)

\end{document}